\documentclass[12pt]{article}

\usepackage{amsmath}
\usepackage{amsfonts}
\usepackage{amssymb}

\newtheorem{theorem}{Theorem}
\newtheorem{corollary}{Corollary}
\newtheorem{definition}{Definition}

\newtheorem{proposition}{Proposition}

\def\F{\mathcal{F}_{q,v}}
\def\R{\mathbb{R}_q^+}

\def\r{\mathbb{R}}
\def\Z{\mathbb{Z}}
\def\N{\mathbb{N}}
\def\I{\infty}
\def\l{\mathcal{L}_{q,1,v}}
\def\L{\mathcal{L}_{q,2,v}}
\def\cb{\mathcal{C}_{q,b}}
\def\M{\mathcal{M}_q^+}
\def\lp{\mathcal{L}_{q,p,v}}
\def\lr{\mathcal{L}_{q,r,v}}
\def\Lp{\mathcal{L}_{q,p',v}}
\def\A{\mathcal{A}_{q,v}}
\def\cc{\mathcal{C}_{q,0}}

\newenvironment{proof}[1][Proof]{\noindent\textbf{#1.} }{\ \rule{0.5em}{0.5em}}

\begin{document}

\title{\bf Functions of $q$-positive type}
\date{ }
\author{Lazhar Dhaouadi \thanks{%
Institut Pr\'eparatoire aux Etudes d'Ing\'enieur de Bizerte
(Universit\'e du 7 novembre à Carthage). Route Menzel Abderrahmene
Bizerte, 7021 Zarzouna, Tunisia. \quad\quad\quad\quad\quad\quad
E-mail lazhardhaouadi@yahoo.fr}} \maketitle

\begin{abstract}
In this paper we characterize the subspace of $\mathcal{L}_{q,1,v}$
of function which are the $q$-Bessel Fourier transform of positive
functions in $\mathcal{L}_{q,1,v}$. As application we give a
$q$-version of the Bochner's theorem.
\end{abstract}

\section{Introduction and Preliminaries}
Given a positive finite Borel measure $\mu$ on the real line $\r$,
the Fourier transform $Q$ of $\mu$ is the continuous function
$$
Q(x)=\int_{\r}e^{-itx}d\mu(t).
$$
The function $Q$ is a positive definite function, i.e for any finite
list of complex numbers $z_1,\ldots,z_n$ and real numbers
$x_1,\ldots,x_n$
$$
\sum_{r=1}^n\sum_{l=1}^n z_r\overline{z_l}Q(x_r-x_l)\geq 0.
$$
Bochner's theorem says the converse is true, i.e. every positive
definite function $Q$ is the Fourier transform of a positive finite
Borel measure. In $q$-Fourier analysis, semelar phenomenon will
appear. It is the subject of our article.

\bigskip

In the following  we consider $0<q<1$ and we adopt the standard
conventional notations of [2]. We put
$$
\mathbb{R}_{q}^{+}=\{q^{n},\quad n\in \mathbb{Z}\},
$$
the set of $q$-real numbers and for complex $a$
$$
(a;q)_{0}=1,\quad (a;q)_{n}=\prod_{i=0}^{n-1}(1-aq^{i}),\quad
n=1...\infty .
$$
Jackson's $q$-integral (see [3]) in the interval $[0,\infty \lbrack
$ is defined by
$$
\int_{0}^{\infty }f(x)d_{q}x=(1-q)\sum_{n=-\infty }^{\infty
}q^{n}f(q^{n}).
$$
Let $\cc$ and $\cb$ denote the spaces of  functions defined on $\R$
continued at $0$, which are respectively vanishing at infinity and
bounded. These spaces are equipped with the topology of uniform
convergence, and by $\lp$ the space of functions $f$ defined on $\R$
such that
$$
 \|f\|_{q,p,v}=\left[\int_0^{\infty}|f(x)|^px^{2v+1}d_qx\right]^{1/p}<\infty.
$$
The $q$-exponential  function is defined by
$$
 e(z,q)=\sum_{n=0}^\infty\frac{z^n}{(q,q)_n}=\frac{1}{(z;q)_\infty},\quad|z|<1.
$$
The normalized Hahn-Exton q-Bessel function of order $v >-1$ (see
[5]) is defined by
$$
j_{v }(z,q)=\sum_{n=0}^{\infty }(-1)^{n}\frac{q^{\frac{n(n-1)}{2}}}{%
(q,q)_{n}(q^{v +1},q)_{n}}z^{n}.
$$
The q-Bessel Fourier transform $\mathcal{F}_{q,v }$ introduced in
[1,4] as follow
$$
\mathcal{F}_{q,v}f(x)=c_{q,v }\int_{0}^{\infty }f(t)j_{v
}(xt,q^{2})t^{2v +1}d_{q}t,
$$
where
$$
c_{q,v }=\frac{1}{1-q}\frac{(q^{2v +2},q^{2})_{\infty }}{(q^{2},q^{2})_{%
\infty }}.
$$
Define the $q-$Bessel translation operator as follows:
$$
T^v_{q,x}f(y)=c_{q,v}\int_0^\infty\mathcal{F}_{q,v}(f)(t)j_v(xt,q^2)j_v(yt,q^2)t^{2v+1}d_qt,\quad\forall
x,y\in\mathbb R_q^+,\forall f\in\mathcal L_{q,v,1}.
$$
Recall that $T^v_{q,x}$ is said positive if $T^v_{q,x}f\geq 0$ for
$f\geq 0$. In the following we tack $q\in Q_v$ where
$$
Q_v=\{q\in ]0,1[,\quad T^v_{q,x}\quad \text{is positive for
all}\quad x\in\R\}.
$$
The $q-$convolution product of both functions $f,g\in\l$ is defined
by
$$
f*_qg(x)=c_{q,v}\int_0^\I T^v_{q,x}f(y)g(y)y^{2v+1}d_qy.
$$
In the end we denote by $\A$ the $q$-Wiener algebra
$$
\A=\left\{f\in\l,\quad \F(f)\in\l\right\}.
$$

\bigskip
The followings results in this sections was proved in [1].
\begin{proposition} Let $n,m\in\Z$ and $n\neq m$, then we
have
$$
c_{q,v}^2 \int_0^\I
j_v(q^nx,q^2)j_v(q^mx,q^2)x^{2v+1}d_qx=\frac{q^{-2n(v+1)}}{1-q}\delta_{nm}.
$$
\end{proposition}

\begin{proposition}
$$
|j_v(q^n,q^2)|\leq\frac{(-q^2;q^2)_\I(-q^{2v+2};q^2)_\I}{(q^{2v+2};q^2)_\I}\left\{
\begin{array}{c}
  1\quad\quad\quad\quad\quad\text{if}\quad n\geq 0 \\
  q^{n^2+(2v+1)n}\quad\text{if}\quad n<0
\end{array}
\right..
$$
\end{proposition}

\begin{proposition}The $q-$Bessel Fourier transform
$$
\F:\l\rightarrow\cc,
$$
satisfying
$$
\|\F(f)\|_{\cc}\leq B_{q,v}\|f\|_{q,1,v},
$$
where
$$
B_{q,v}=\frac{1}{1-q}\frac{(-q^2;q^2)_\I(-q^{2v+2};q^2)_\I}{(q^2;q^2)_\I}.
$$
\end{proposition}

\begin{theorem}
 Given $f\in\l$ then we have
$$
\F^2(f)(x)=f(x),\quad\forall x\in\R.
$$
If $f\in\l$ and $\F(f)\in\l$ then
$$
\|\F(f)\|_{q,2,v}=\|f\|_{q,2,v}.
$$
\end{theorem}

\begin{proposition} Let $f\in\l$ then
$$
T^v_{q,x}f(y)=\int_0^\I f(z)D_v(x,y,z)z^{2v+1}d_qz,
$$
where
$$
D_v(x,y,z)=c^2_{q,v}\int_0^\I
j_v(xt,q^2)j_v(yt,q^2)j_v(zt,q^2)t^{2v+1}d_qt.
$$
\end{proposition}

\begin{proposition} Given two functions $f,g\in\mathcal L_{q,v,1}$ then
$$
f*_qg\in\mathcal L_{q,v,1},
$$
and
$$
\mathcal{F}_{q,v}(f*_qg)=\mathcal{F}_{q,v}(f)\times\mathcal{F}_{q,v}(g).
$$
\end{proposition}

\begin{proposition}The $q$-Gauss kernel
$$
G^v(x,t,q^2)=\frac{(-q^{2v+2}t,-q^{-2v}/t;q^2)_\infty}{(-t,-q^2/t;q^2)_\infty}e(-\frac{q^{-2v}}{t}x^2,q^2),
$$
satisfying
$$
\mathcal F_{q,v}\left\{e(-ty^2,q^2)\right\}(x)=G^v(x,t,q^2),
$$
and for all function $f\in\cb$
$$
\lim_{a\rightarrow
0}c_{q,v}\int_0^{\infty}f(x)G^v(x,a^2,q^2)x^{2v+1}d_qx=f(0).
$$
\end{proposition}

\begin{theorem} Given $1<p,p',r\leq 2$ and
$$
{1\over p}+{1\over p'}-1={1\over r}.
$$
If $f\in\lp$ and $g\in\Lp$ then
$$
f*_qg\in\lr.
$$
\end{theorem}

\section{Functions of $q$-positive type}

\begin{definition}A function $\phi$ is of $q$-positive type if
$$
\phi\in\cb\cap\l
$$
and for any finite list of complex numbers $z_1,\ldots,z_n$ and
$q$-real numbers $x_1,\ldots,x_n$
\begin{equation}
\sum_{r=1}^n\sum_{l=1}^n z_r\overline{z_l}T_{q,x_r}^v\phi(x_l)\geq
0.
\end{equation}
\end{definition}

\begin{proposition}Let $\phi\in\A$ of $q$-positive type then $\F\phi$  is of $q$-positive type.
\end{proposition}

\begin{proof}From Proposition 3 and the definition of the
$q$-Wiener algebra
$$
\F(\phi)\in\cb\cap\l.
$$
On the other hand, with the inversion formula in Theorem 1 we get
$$
T^v_{q,x}\F(\xi)(y)=\int_0^{\infty}j_v(tx,q^2)j_v(ty,q^2)t^{2v+1}\phi(t)d_qt,
$$
then
$$\aligned
\sum_{r=1}^n\sum_{l=1}^n
z_r\overline{z_l}T_{q,x_r}^v\F\xi(x_l)&=c_{q,v}\int_0^\I
\left[\sum_{r=1}^n\sum_{l=1}^n z_r\overline{z_l}
j_v(x_rt,q^2)j_v(x_lt,q^2)\right]t^{2v+1}\phi(t)d_qt\\
&=c_{q,v}\int_0^\I\left[\sum_{r=1}^n z_r
j_v(x_rt,q^2)\right]\overline{\left[\sum_{l=1}^n z_l
j_v(x_lt,q^2)\right]}t^{2v+1}\phi(t)d_qt\\
&=c_{q,v}\int_0^\I\left|\sum_{r=1}^n z_r
j_v(x_rt,q^2)\right|^2t^{2v+1}\phi(t)d_qt\geq 0.
\endaligned$$
This finish the proof.
\end{proof}

\begin{proposition}If $\phi$ is of $q$-positive type and $f\in\L$ then
$$
\phi*_qf\in\L,
$$
and
$$
\langle \phi*_qf,f\rangle\geq 0.
$$
\end{proposition}

\begin{proof}From Theorem 2 we see that $\phi*_qf\in\L$. On the
other hand
$$\aligned
\langle \phi*_qf,f\rangle &=c^2_{q,v}\int_0^\I\left[\int_0^\I
T^v_{q,x}\phi(y)f(y)y^{2v+1}d_qy\right]f(x)x^{2v+1}d_qx\\
&=(1-q)^2c^2_{q,v}\sum_{r=1}^\I\sum_{l=1}^\I
q^{(2v+2)r}f(q^r)q^{(2v+2)l}f(q^l)T^v_{q,q^r}\phi(q^l)\geq 0.
\endaligned$$
This finish the proof.
\end{proof}

\begin{corollary}If $\phi$ is of $q$-positive type then
$$
\F\phi(x)\geq 0,\quad\forall x\in\R.
$$
\end{corollary}

\begin{proof}Given $x\in\R$ and let
$$
f_x:t\mapsto c_{q,v}j_v(xt,q^2),
$$
then with Proposition 2 we see that $f_x\in\L$ and by Proposition 1
$$
\F f_x(y)=\delta_{q,v}(x,y),
$$
which implies (see[1])
$$
\langle \F\phi\times\F f_x,\F
f_x\rangle=\F\phi(x)\delta_{q,v}(x,x)=\frac{1}{(1-q)x^{2v+2}}\F\phi(x).
$$
From Proposition 8
$$
\langle \F\phi\times\F f,\F f\rangle=\langle \phi*_qf,f\rangle\geq
0,
$$
this leads to the result.
\end{proof}

\begin{proposition}If $\phi$ is of $q$-positive type then
$\F\phi\in\l$.
\end{proposition}

\begin{proof}
From Proposition 6
$$
\lim_{a\rightarrow 0}\int_0^\I
e(-a^2x,q^2)\F\phi(x)x^{2v+1}d_qx=\lim_{a\rightarrow
0}c_{q,v}\int_0^\I G^v(x,a^2,q^2)\phi(x)x^{2v+1}d_qx=\phi(0).
$$
By the monotone convergence theorem and the preview corollary we see
that
$$
\int_0^\I |\F\phi(x)|x^{2v+1}d_qx=\int_0^\I
\F\phi(x)x^{2v+1}d_qx=\phi(0).
$$
This finish the proof.
\end{proof}

\begin{corollary}If $\phi$ is of $q$-positive type then there exist a positive function $\xi\in\A$ such that
$$
\phi(x)=\F\xi(x),\quad\forall x\in\R.
$$
\end{corollary}

\begin{proof}From the inversion formula in theorem 1
$$
\phi(x)=\F^2\phi(x),\quad\forall x\in\R.
$$
Define the function $\xi$ as follows
$$
\xi(x)=\F\phi(x).
$$
By the use of Corollary 1 and Proposition 9 we see that $\xi$ is a
positive function of $\A$.
\end{proof}

\begin{proposition}Suppose $\phi$ is of $q$-positive type.
If $f\in\l$ is positive function then the product $\phi\F f$ is of
$q$-positive type.
\end{proposition}

\begin{proof}Proposition 5 and Proposition 9 give
$$
\F(\phi\F f)(t)=\F\phi*_q f(t),\quad\forall t\in\R,
$$
then
$$\aligned
&\sum_{r=1}^n\sum_{l=1}^n z_r\overline{z_l}T_{q,x_r}^v\left(\phi\F
f\right)(x_l)\\
&=c_{q,v}\int_0^\I\left[\sum_{r=1}^n\sum_{l=1}^n
z_rj_v(x_rt,q^2)\overline{z_lj_v(x_lt,q^2)}\right]\F(\phi\F f)(t)t^{2v+1}d_qt\\
&=c_{q,v}\int_0^\I\left[\sum_{r=1}^n\sum_{l=1}^n
z_rj_v(x_rt,q^2)\overline{z_lj_v(x_lt,q^2)}\right]\F\phi*_q f(t)t^{2v+1}d_qt\\
&=c_{q,v}\int_0^\I\left|\sum_{r=1}^n
z_rj_v(x_rt,q^2)\right|^2\F\phi*_q f(t)t^{2v+1}d_qt.
\endaligned$$
From the definition of the $q$-convolution product we write
$$
\F\phi*_q f(t)=c_{q,v}\int_0^\I \F\phi(z) T_{q,t}f(z)z^{2v+1}d_qz.
$$
Proposition 4 give
$$
T_{q,t}f(z)=c_{q,v}\int_0^\I D_v(t,z,s)f(s)s^{2v+1}d_qs\geq 0.
$$
This implies with Corollary 1
$$
\F\phi*_q f(t)\geq 0,
$$
which leads to the result.
\end{proof}

\begin{corollary}Given two functions $\phi_1,\phi_2$ which are of $q$-positive type
 then the product $\phi_1\times\phi_2$ is also of $q$-positive type.
\end{corollary}

\begin{proof}Let $\xi=\F\phi_2$ then with the inversion formula in
theorem 1 we see that $\F\xi=\phi_2$. Proposition 9 give
$$
\xi\in\l,
$$
and by Proposition 10 we achieved the proof.
\end{proof}

\section{$q$-Bochner's Theorem}
We consider the set $\M$ of positives and bonded measures on $\R$.
The $q$-Bessel Fourier transform of $\xi\in\M$ is defined by
$$
\F(\xi)(x)=\int_0^{\infty}j_v(tx,q^2)t^{2v+1}d_q\xi(t).
$$
The $q-$convolution product of two measures $\xi,\rho\in\M$ is given
by
$$
\xi*_q\rho(f)=\int_0^{\infty}T^v_{q,x}f(t)t^{2v+1}d_q\xi(x)d_q\rho(t),
$$
and we have
$$
\F(\xi*_q\rho)=\F(\xi)\F(\rho).
$$
The following Theorem (see[1]) is crucial for the proof of our main
result.
\begin{theorem} Let $(\xi_n)_{n\geq 0}$ be a sequences of probability measures of
$\M$ such that
$$
\lim_{n\rightarrow\infty}\F(\xi_n)(x)=\psi(x),
$$
then there exists $\xi\in\M$ such that the sequence $\xi_n$ converge
strongly toward $\xi$ and
$$
\F(\xi)=\psi.
$$
\end{theorem}

In the following we consider the function $\psi$ defined by
$$
\psi(x)=\left\{\begin{array}{c}
          1-x\quad\text{if}\quad x<1 \\
          0\quad\text{otherwise}
        \end{array}\right..
$$
Now we are in a position to state and prove the q-analogue of the
Bochner's theorem

\begin{theorem}Let $\phi$ be a function defined on $\R$ continued at
$0$. Assume that the following function
$$
\phi_n:x\mapsto\phi(x)\psi(q^nx),
$$
satisfy (1) for all $n\in\N$ then there exist $\xi\in\M$ such that
$$
\F(\xi)=\phi.
$$
\end{theorem}

\begin{proof}The function $\phi_n$ is of $q$-positive type. From
Corollary 2 there exist $\varrho_n$ a positive function of $\A$ such
that
$$
\F(\varrho_n)=\phi_n.
$$
The measure $\xi_n$ defined by
$$
d_q\xi_n(x)=\varrho_n(x)d_qx,
$$
belong to $\M$ and
$$
\int_0^\I x^{2v+1}d_q\xi_n(x)=\F(\varrho_n)(0)=\phi_n(0)=\phi(0).
$$
Assume that $\phi(0)=1$. On the other hand
$$
\lim_{n\rightarrow\infty}\F(\xi_n)(x)=\lim_{n\rightarrow\infty}\phi_n(x)=\phi(x).
$$
From Theorem 3 there exists $\xi\in\M$ such that the sequence
$\xi_n$ converge strongly toward $\xi$, and
$$
\F(\xi)=\phi,
$$
which leads to the result.
\end{proof}

\end{document}